\documentclass[12pt]{amsart}
\usepackage{amssymb,mathrsfs,fullpage}

\input xy
\xyoption{all}

\newcommand{\g}{\mathfrak g}

\renewcommand{\t}{\mathfrak t}

\renewcommand{\l}{\mathfrak l}

\renewcommand{\u}{\mathfrak u}
\newcommand{\p}{\mathfrak p}

\newcommand{\so}{\mathfrak{so}}
\renewcommand{\sp}{\mathfrak{sp}}

\newcommand{\sub}{\subseteq}

\renewcommand{\bar}{\overline}

\newcommand{\C}{\mathbb C}

\newcommand{\Z}{\mathbb Z}

\newcommand{\cO}{\mathcal O}

\DeclareMathOperator{\GL}{GL} %
\DeclareMathOperator{\SL}{SL} %
\DeclareMathOperator{\Sp}{Sp} %
\DeclareMathOperator{\SO}{SO} %
\DeclareMathOperator{\OR}{O} %

\DeclareMathOperator{\Char}{char} %
\DeclareMathOperator{\Lie}{Lie} %
\DeclareMathOperator{\rank}{rank} %

\numberwithin{equation}{section}

\newtheorem{thm}{Theorem}[section]

\theoremstyle{definition}
\newtheorem{defn}[thm]{Definition}
\newtheorem{exmp}[thm]{Examples}

\theoremstyle{remark}
\newtheorem{rem}[thm]{Remark}

\title[Richardson elements for classical groups]
{\boldmath Richardson elements for parabolic subgroups of classical
groups in positive characteristic}
\author[K.~Baur and S.~M.~Goodwin]
{Karin Baur and Simon M.~Goodwin}
\date{30.08.2006}

\address{Department of Mathematics,
University of Leicester, Leicester, LE1 7RH, UK}
\email{k.baur@mcs.le.ac.uk}
\urladdr{http://www.math.le.ac.uk/people/kbaur}

\address{School of Mathematics, University of Birmingham, Birmingham, B15 2TT, UK}
\email{goodwin@maths.bham.ac.uk}
\urladdr{http://web.mat.bham.ac.uk/S.M.Goodwin}


\begin{document}

\begin{abstract}
Let $G$ be a simple algebraic group of classical type over an
algebraically closed field $k$. Let $P$ be a parabolic subgroup of
$G$ and let $\p = \Lie P$ be the Lie algebra of $P$ with Levi
decomposition $\p = \l \oplus \u$, where $\u$ is the Lie algebra of
the unipotent radical of $P$ and $\l$ is a Levi complement. Thanks
to a fundamental theorem of R.~W.~Richardson (\cite{Ri}), $P$ acts
on $\u$ with an open dense orbit; this orbit is called the {\em
Richardson orbit} and its elements are called {\em Richardson
elements}. Recently (\cite{Ba}), the first author gave constructions
of Richardson elements in the case $k = \C$ for many parabolic
subgroups $P$ of $G$. In this note, we observe that
these constructions remain valid for any algebraically closed field
$k$ of characteristic not equal to 2 and we give constructions of
Richardson elements for the remaining parabolic subgroups.
\end{abstract}

\maketitle

\section{Introduction}

Let $G$ be a reductive algebraic group over an algebraically closed
field $k$ and let $P$ be a parabolic subgroup of $G$ with unipotent
radical $U$.  We write $\g$, $\p$ and $\u$ for the Lie algebras of
$G$, $P$ and $U$ respectively.  It is well-known that $G$ has
finitely many nilpotent orbits in $\g$: this was first proved by
R.~W.~Richardson (\cite{Ri1}), when $\Char k$ is zero or good for
$G$ (see also \cite[Thm.\ 5.4]{SS}); we refer the reader to
\cite[5.11]{Ca} for a survey of the result in bad characteristic. It
follows that there is a unique nilpotent orbit $G \cdot e$ which
intersects $\u$ in an open dense subvariety. Richardson's dense
orbit theorem (\cite{Ri}) tells us that the intersection $G \cdot e
\cap \u = P \cdot e$ is a single $P$-orbit (we may assume that $e
\in \u$). The $P$-orbit $P \cdot e$ is called the {\em Richardson
orbit} and its elements are called {\em Richardson elements}; note
that we also sometimes refer to $G \cdot e$ as the Richardson orbit
of $P$.

The purpose of this paper is to give explicit constructions
of Richardson elements $e \in \u$ in case $G$ is of classical type
and $\Char k \ne 2$; this extends work of the first author in
\cite{Ba}. Further, the representatives that we construct are of a
``nice'' form as explained in Theorem \ref{T:support} below.  We
only consider orthogonal and symplectic groups in this paper,
because a construction of Richardson elements for $G$ a special
linear (or general linear) group was given in \cite{BHRR}.

We note that Richardson's dense orbit theorem also applies to the
action of $P$ on its unipotent radical $U$ by conjugation, i.e.\
there is an open dense $P$-conjugacy class in $U$.  If $\Char k$ is
good for $G$ then, under minor restrictions on $G$, there exists a
Springer isomorphism, i.e.\ a $G$-equivariant isomorphism from the
unipotent variety of $G$ to the nilpotent variety of $\g$, see
\cite[III, 3.12]{SS} and \cite[Cor.\ 9.3.4]{BR}.  A Springer
isomorphism restricts to a $P$-equivariant isomorphism of varieties
$\phi : U \to \u$; in fact such a $P$-equivariant isomorphism exists
without any restrictions on $G$. Therefore, our constructions of
Richardson elements $e \in \u$ can be used to give Richardson
elements in $U$. In case $G$ is an orthogonal or symplectic group,
we can take a Cayley map for a Springer isomorphism, see \cite[III,
3.14]{SS}.

\medskip

In case $G = \GL_n(k)$ or $\SL_n(k)$ a construction of Richardson
elements was given by T.~Br\"ustle, L.~Hille, C.~M.~Ringel and
G.~R\"ohrle in \cite[\S 8]{BHRR} valid for any algebraically closed
field $k$.  This construction was a consequence of results
concerning the representation theory of a certain quasi-hereditary
algebra and defines the Richardson elements from certain {\em line
diagrams} in the plane.  It is also recalled in \cite[\S 3]{Ba} and
a direct proof that the construction yields Richardson elements is
given.

In \cite{Ba} the first author constructed Richardson elements for
many parabolic subgroups of a classical group over $k = \C$ using
line diagrams generalizing those in \cite{BHRR}.  The purpose of
this paper is to extend these constructions to any parabolic
subgroup of a classical group and further to justify that the
constructions remain valid for any algebraically closed field $k$ of
characteristic not equal to 2.  We note that the line diagrams
defined in this paper are similar to the {\em pyramids} defined in
\cite[\S 5 and 6]{EK}.

In case $G$ is an exceptional group, R.~Lawther has given
representatives of all Richardson elements $u \in U$ in the
unpublished article \cite{La}; one can check that the tables in {\em
loc.\ cit.\ } also determine Richardson elements in $x \in \u$, when
$\Char k$ is good for $G$. A.~G.~Elashvili has also constructed
representatives; these constructions were made when Elashvili
determined the Lusztig--Spaltenstein induction (see \cite{LS}) for
nilpotent orbits in exceptional groups, this work is reported in the
appendix to Chapter 2 of \cite{Sp}. In addition, in \cite{Ba1} the
first author gave representatives of Richardson elements for {\em
nice} parabolic subalgebras.

\medskip

As a consequence of our constructions in this paper along with those
in \cite{BHRR} and \cite{Ba}, and the representatives given in
\cite{La}, we have Theorem \ref{T:support} below, which says that
one can always find Richardson elements of a ``nice'' form. Before
stating this theorem we need to introduce some notation.

Let $G$ be a simple algebraic group over an algebraically closed
field $k$.  Assume that $\Char k$ is zero or a good prime for $G$.
Let $T$ be maximal torus of $G$ and $P$ a parabolic subgroup of $G$
containing $T$. Let $\Phi$ be the root system of $G$ with respect to
$T$. For $\alpha \in \Phi$, let $\g_\alpha$ denote the root subspace
of $\g = \Lie G$ corresponding to $\alpha$ and let $e_\alpha$ be a
generator of $\g_\alpha$. Let $\u$ be the Lie algebra of the
unipotent radical of $P$ and let $\Phi(\u) \sub \Phi$ be defined by
$\u = \bigoplus_{\alpha \in \Phi(\u)} \g_\alpha$.  Given $x \in \g$,
we write $C_G(x)$ for the centralizer of $x$ in $G$.  For a closed
subgroup $H$ of $G$, we write $N_G(H)$ and $C_G(H)$ for the normalizer
and centralizer of $H$ in $G$ respectively.

We can now state Theorem \ref{T:support}.  It is straightforward to
verify the first assertion for classical groups from the
constructions in \cite{BHRR}, \cite{Ba} and this paper; and from the
tables in \cite{La} for exceptional groups.  We discuss the claim
about the minimality of the size of $\Gamma$ for $G$ an orthogonal
or symplectic group below; it is undemanding to verify the claim for
$G$ a special linear or exceptional group from the representatives
given in \cite{BHRR} and \cite{La}.

\begin{thm} \label{T:support}
There exists a subset $\Gamma$ of $\Phi(\u)$ consisting of linearly
independent roots such that
$$
x = \sum_{\alpha \in \Gamma} e_\alpha \in \u
$$
is a Richardson element for $P$.  Moreover, we can find $\Gamma$
with $|\Gamma| = \rank G - \rank C_G(x)$ and this is the minimal
possible size of $\Gamma$.
\end{thm}

We call a representative of the Richardson orbit for $P$ of the form
$x = \sum_{\alpha \in \Gamma} a_\alpha e_\alpha \in \u$, where
$a_\alpha \in k^\times$, and $\Gamma$ is a set of linearly
independent roots with $|\Gamma| = \rank G - \rank C_G(x)$, a {\em
minimal Richardson element}. We do not make any assertion about
uniqueness of minimal Richardson elements here. Indeed let $L$ be
the Levi subgroup of $P$ containing $T$,
then for a minimal Richardson element $x$ and $g \in N_L(T)$, we
have that $g \cdot x$ is a minimal Richardson element.  In case $x$
is regular nilpotent in the Lie algebra of the Levi subgroup
$C_G(S)$ of $G$, where $S$ is a maximal torus of $C_G(x)$, one can
in fact prove that all minimal Richardson elements for $P$ are of
the form $g \cdot x$ for some $g \in N_L(T)$.  This need not be the
case in general.

It is unclear whether Theorem \ref{T:support} holds without the
assumption that $\Char k$ is good for $G$.  We note that for
orthogonal groups the constructions given in the present paper, need
not yield Richardson elements when $\Char k = 2$.  For the
symplectic groups the representatives for Richardson orbits that we
give do remain valid for $\Char k = 2$; though this does not follow
from the proofs we give in this article.


Let $S$ be as in the statement of Theorem \ref{T:support} and
further assume that $S \sub T$.  According to the Bala--Carter
classification of nilpotent orbits (see \cite[\S 4]{Ja}), $x$ is
distinguished nilpotent in the Levi subgroup $C_G(S)$ of $G$.   In
particular, this forces $\Gamma$ to have size at least the
semisimple rank of $C_G(S)$, i.e.\ $\rank T - \rank S$.  To verify
that the Richardson elements $x$ given in this paper are minimal we
determine the type of the orbit of $x$. We recall that the {\em
type} of (the orbit of) $x$ is by definition the conjugacy class in
$G$ of the derived subgroup $C_G(S)'$ of $C_G(S)$; this is usually
determined by the isomorphism type of $C_G(S)'$.  For the types that
we give, it is rarely the case that they are not determined by just
the isomorphism type of $C_G(S)'$, we make brief remarks where there
are different conjugacy classes.  We check that for our
constructions we have that $|\Gamma|$ is equal to the rank of
$C_G(S)'$; this is carried out at the end of subsections
\ref{ss:orthgen} and \ref{ss:sympgen}.  The type of $x$ forms part
of the Bala--Carter label of the nilpotent orbit of $x$.  We note
that it is straightforward to determine the Bala--Carter labels for
the Richardson orbits, but this is technical so we choose not to
include it in this paper.

\medskip

The minimality of the size of $\Gamma$ for the representatives
constructed in this paper along with those given in \cite{BHRR} and
\cite{La} allows one to give an alternative proof of the
classification of {\em nice parabolic subalgebras} given by
N.~Wallach and the first author in \cite{BW}. Further, one can show
that the classification of nice parabolic subalgebras given in {\em
loc.\ cit.\ }remains valid in good positive characteristic. We
explain how this can be achieved below.

Let $\g = \sum_{i \in \Z} \g_i$ be the $\Z$-grading of $\g$
associated to $\p$ (and $\t$), see for example \cite[\S 1]{BW}. So
we have that $\p = \sum_{i \ge 0} \g_i$, $\u = \sum_{i \ge 1} \g_i$.
We recall that $\p$ is called {\em nice} if there exists a
Richardson element $x \in \g_1$.

Let $x = \sum_{\alpha \in \Gamma} e_\alpha \in \u$ be a minimal
Richardson element for $P$. Let $\Gamma_1$ be the subset of $\Gamma$
consisting of roots $\alpha$ such that $\g_\alpha \sub \g_1$ and let
$x_1 = \sum_{\alpha \in \Gamma_1} e_\alpha$.  It follows from a
result of Richardson (\cite[Thm.\ E]{Ri3}) that $L$ acts on $\g_1$
with a dense orbit, where $L$ is the Levi subgroup of $P$ containing
$T$.  It follows from the fact that $P \cdot x$ is dense in $\u$
that $L \cdot x_1$ is dense in $\g_1$. Suppose that $\p$ is nice and
let $y \in \g_1$ be a Richardson element. Then $L \cdot y$ is dense
in $\g_1$ and it follows that $y$ is in the same $P$-orbit as $x_1$,
thus $x_1$ is Richardson.  The minimality of $|\Gamma|$ implies that
$\Gamma_1 = \Gamma$.  Hence, we deduce that $\p$ is a nice parabolic
subalgebra if and only if $\g_\alpha \sub \g_1$ for all $\alpha \in
\Gamma$. From this equivalence and the representatives of Richardson
orbits given in \cite{BHRR}, \cite{La} and this paper, one can
determine which parabolic subalgebra are nice; and that it does not
depend on the characteristic of $k$, for $\Char k$ zero or good for
$G$.

\medskip

In the language of \cite{H}, a nilpotent $G$-orbit is {\em
polarizable} if it is the Richardson orbit for some parabolic
subgroup of $G$. As a consequence of our constructions one can
easily obtain the classification of polarizable nilpotent $G$-orbits
for classical $G$ given in Theorem \ref{T:orbits} below.  This
classification is also given by G.~Kempken in \cite[\S 4]{K} as a
consequence of the determination of the Lusztig--Spaltenstein
induction for nilpotent orbits in classical groups presented in
\cite{K}.  The Lusztig--Spaltenstein induction for nilpotent orbits
in classical groups is also described in \cite[Ch.\ 2, \S 7]{Sp}. In
the case $G = \SL_n(k)$ it is well-known that every nilpotent orbit
is the Richardson orbit of some parabolic subgroup. Further, we note
that one can obtain the classification for exceptional groups from
the tables in \cite{La}.

We require some notation for the statement of Theorem
\ref{T:orbits}. We recall that nilpotent orbits for classical groups
are given by certain partitions (given by Jordan normal forms).  A
partition is a sequence $\lambda =
(\lambda_1,\lambda_2,\dots,\lambda_l)$ of non-increasing positive
integers; $\lambda$ is said to be a partition of $\sum_{i=1}^l
\lambda_i$. The nilpotent orbits of the orthogonal group $\OR_N(k)$
are parametrized by {\em orthogonal partitions} of $N$, i.e.\
partitions of $N$ for which all even parts occur with even
multiplicity; the nilpotent orbits of the symplectic group
$\Sp_{2n}(k)$ are parametrized by {\em symplectic partitions} of
$2n$, i.e.\ partitions of $2n$ for which all odd parts occur with
even multiplicity.

Given partitions $\lambda^i = (\lambda^i_1,\dots,\lambda^i_{l_i})$
(for $i = 1,\dots,m$) with $\lambda^i_{l_i} \ge \lambda^{i+1}_1$ for
each $i$, we write $(\lambda^1,\dots,\lambda^m)$ for the partition
$(\lambda^1_1,\dots,\lambda^1_{l_1},\lambda^2_1,\dots,\lambda^m_{l_m})$.
Let $\lambda = (\lambda_1,\dots,\lambda_l)$ be a partition and
suppose that $l$ is even. We say that $\lambda$ is: an {\em
even-pair partition} if each $\lambda_i$ is even and $\lambda_{2i-1}
= \lambda_{2i}$ for each $i$; an {\em odd-pair partition} if each
$\lambda_i$ is odd and $\lambda_{2i-1} = \lambda_{2i}$ for each $i
> 0$; an {\em even 2-step-descending partition} if each $\lambda_i$ is even
and $\lambda_{2i} > \lambda_{2i+1}$ for each $i$; an {\em odd
2-step-descending partition} if each $\lambda_i$ is odd and
$\lambda_{2i} > \lambda_{2i+1}$ for each $i$.

\begin{thm} \label{T:orbits}
Assume $\Char k \ne 2$.
\begin{enumerate}
\item [(i)] Let $G = \OR_N(k)$ and let $\cO_\lambda$ be the nilpotent
orbit given by the orthogonal partition $\lambda$.  Then
$\cO_\lambda$ is the Richardson orbit for some parabolic subgroup of
$G$ if and only if $\lambda$ is of the form
$$
\lambda = (\lambda^0,\lambda^1,\dots,\lambda^m),
$$
where: $m \ge 0$; $\lambda^0$ is a (possibly empty) partition with
all odd entries; for odd $i > 0$, $\lambda^i$ is an even-pair
partition; and for even $i
> 0$, $\lambda^i$ is an odd 2-step-descending partition.
\item[(ii)] Let $G = \Sp_{2n}(k)$ and let
$\cO_\lambda$ be the nilpotent orbit given by the symplectic
partition $\lambda$. Then $\cO_\lambda$ is the Richardson orbit for
some parabolic subgroup of $G$ if and only if $\lambda$ is of the
form
$$
\lambda = (\lambda^1,\dots,\lambda^m,\mu),
$$
where: $m \ge 0$; for odd $i
> 0$, $\lambda^i$ is an odd-pair partition ($\lambda^1$ is allowed to be the empty partition); for even $i > 0$,
$\lambda^i$ is an even 2-descending partition; and $\mu$ is a
(possibly empty) partition with all even entries.
\end{enumerate}
\end{thm}

We now outline the structure of this paper.  In Section
\ref{S:prelims} we give the required recollection about nilpotent
orbits in classical groups.
In Section \ref{S:orth} we give the construction of Richardson
elements for parabolic subgroups of $\OR_N(k)$. Finally, in Section
\ref{S:symp} we give the construction of Richardson elements for
parabolic subgroups of $\Sp_{2n}(k)$.  We give complete proofs that
the construction is correct for the $\OR_N(k)$ case; the
$\Sp_{2n}(k)$ case is similar (and easier) so we omit the proofs in
this case.

\section*{Acknowledgments}

The first author was supported by the EPSRC grant number
GR/S35387/01. The second author thanks New College, Oxford for
financial support whilst this work was completed.  We are grateful
to G.~R\"ohrle for reading through an earlier version and providing
corrections.  We thank A.~G.~Elashvili and R.~I.~Lawther for
comments about Theorem \ref{T:support}.

\section{Preliminaries} \label{S:prelims}

Throughout this paper $k$ is an algebraically closed field of
characteristic not equal to 2.

By a partition $\lambda$ we mean a sequence $\lambda =
(\lambda_1,\dots,\lambda_l)$ of non-increasing positive integers; we
say $\lambda$ is a partition of $\sum_{i=1}^l \lambda_i$.  We recall
the dominance ordering on partitions is defined by $\lambda \le \mu$
if and only if $\sum_{i=1}^j \lambda_i \le \sum_{i=1}^j \mu_i$ for
all $j$.

We recall that the nilpotent orbits of the general linear group $G =
\GL_n(k)$ are parametrized by the partitions of $n$.  Let $x$ and
$y$ lie in the nilpotent $G$-orbits parametrized by $\lambda$ and
$\mu$ respectively.  Then it is well-known that $G \cdot x \sub
\overline{G \cdot y}$ if and only $\lambda \le \mu$, see \cite[\S
3]{H2}.

\subsection{Orthogonal groups}
\label{ss:orth}  Let $G = \OR_N(k)$ and $\g = \Lie G = \so_N(k)$. We
assume $N \ge 3$ and let $n = \lfloor \frac{N}{2}\rfloor$. Let $V =
k^N$ be the natural $N$-dimensional $G$-module with standard basis
$v_{-n},\dots,v_{-1},v_0,v_1,\dots,v_n$ and $G$-invariant symmetric
bilinear form $(\,,)$ defined by $(v_0, v_i) = (v_0, v_{-i}) = 0,
(v_0,v_0) = 1$, $(v_i, v_j) = (v_{-i},v_{-j}) = 0$ and $(v_i,
v_{-j}) = \delta_{i,j}$ for $1 \leq i,j \leq n$ (omitting $v_0$
everywhere if $N$ is even).  Let $T$ be the maximal torus of $G$,
which acts diagonally on the standard basis of $V$ and let $\Phi$ be
the root system of $G$ with respect to $T$.  The following matrices
give a Chevalley basis for $\g = \so_N(k)$ (omitting the last family
if $N$ is even):
\begin{align*}
\{e_{i,j} - e_{-j,-i}\}_{1 \leq i,j \leq n} &\cup
\{e_{i,-j}-e_{j,-i}, e_{-j,i}-e_{-i,j}\}_{1 \leq i < j \leq n}\\
&\cup \{e_{k,0}-e_{0,-k}, e_{0,k}-e_{-k,0}\}_{1 \leq k \leq n},
\end{align*}
where $e_{i,j}$ denotes the $ij$-matrix unit for $-n \le i,j \le n$,
i.e.\ $e_{i,j}$ is the linear map of $V$ sending $v_i$ to $v_j$ and
$v_{i'}$ to $0$ for $i' \ne i$. We define $\epsilon_{i,j} = \pm 1$
to be the coefficient of $e_{i,j}$ in this Chevalley basis if it
appears.

Let $d' = (d_1,\dots,d_s) \in \Z_{\ge 0}^s$ satisfy $\sum_{i=1}^s
d_i \le n$. Let $t = 2s+1$, and define $d =
(d_{-s},d_{-s-1},\dots,d_{s-1},d_s) \in \Z_{\ge 0}^t$ by setting $d_0 = N -
2\sum_{i=1}^s d_i$, and $d_{-i} = d_i$ for $i = 1,\dots,s$.  A
vector $d$ as constructed above is called an {\em orthogonal
dimension vector}; we say that $d$ is {\em proper} if: $d_i \ne 0$
for all $i \ne 0$; and if $2(\sum_{i=1}^s d_i) = N$, then $d_1 \ne
1$.  We write $|d| = \sum_{i=-s}^s d_i$.

Let $d$ be an orthogonal dimension vector.  We define $c_i =
\sum_{j=0}^{i-1} d_j$ for $i = 1,\dots,s$. For $i = 1,\dots,s+1$,
let $V_{-i}$ be the subspace of $V$ generated by
$v_{-n},v_{-n+1},\dots,v_{-c_i-1}$. The stabilizer $P(d)$ of the
isotropic flag $0 = V_{-(s+1)} \sub V_{-s} \sub \dots \sub V_{-1}$
is a parabolic subgroup of $G$. Moreover, for any parabolic subgroup
$P$ of $G$, there is a unique proper orthogonal dimension vector
$d$ such that $P$ is conjugate
to $P(d)$.  Therefore, in what follows it suffices to just consider
parabolic subgroups of the form $P(d)$. We note that if $d$ is
proper, then it gives the sizes of the blocks in the Levi subgroup
of $P$ containing $T$.

If $\tilde d$ differs from $d$ by removing zero entries, then we
have $P(\tilde d) = P(d)$. Further, if we have $\tilde d_1 = \tilde
d_{-1} = 1$ and $\tilde d_0 = 0$ and $d$ is obtained by replacing
the 3-tuple $(1,0,1)$ in the centre of $\tilde d$ with the single
entry 2, then $P(\tilde d) = P(d)$. This explains why all parabolic
subgroups are obtained  by just considering proper $d$.  We allow
$d$ not to be proper as we need to consider such $d$ in the
inductive construction of Richardson elements given in Definition
\ref{D:orth}.

We may identify $\GL_N(k)$ with $\GL(V)$. For $i > 0$, we define
$V_i$ to be the subspace generated by
$v_{-n},\dots,v_{-1},v_0,v_1,\dots,v_{c_i}$ (we omit $v_0$ if $N$ is
even). The stabilizer $Q(d)$ of the flag $0 = V_{-(s+1)} \sub
V_{-s+1} \sub \dots \sub V_{-1} \sub V_1 \sub \dots \sub V_{s+1}$ in
$GL_N(k)$ is a parabolic subgroup of $\GL_N(k)$ and moreover we have
$P(d) = Q(d) \cap G$.

We recall that the nilpotent $G$-orbits are parametrized by
partitions $\lambda$ of $N$ such that every even part of $\lambda$
appears with even multiplicity. We call such a partition an {\em
orthogonal partition}.  Let $x$ and $y$ lie in the nilpotent
$G$-orbits parametrized by the orthogonal partitions $\lambda$ and
$\mu$ respectively. Then $G \cdot x \sub \overline{G \cdot y}$ if
and only $\lambda \le \mu$, see \cite[\S 3]{H2}.  This means that
for nilpotent $x,y \in \g$ we have $x \in \overline{G \cdot y}$ if
and only if $x \in \overline{\GL_N(k) \cdot y}$.

\begin{rem} \label{R:OnotSO}
We note that if we considered the group $G = \SO_N(k)$ instead of
$\OR_N(k)$, then the description of the conjugacy classes of
parabolic subgroups of $G$ and nilpotent $G$-orbits is slightly more
complicated. In case $N$ is even, some conjugacy classes of
parabolic subgroups of $\OR_N(k)$ give rise to two conjugacy classes
of parabolic subgroups of $G$; and if $N$ is divisible by 4, then some
nilpotent orbits for $\OR_N(k)$ split into two $G$-orbits.

Let $P = P(d)$ be a parabolic subgroup of $\OR_N(k)$.  Then either
$P \sub G$ or $P \cap G$ is a subgroup of $P$ of index 2. Let $x$ be
a representative of the dense $P$-orbit in $\u$, then we must have
$P \cdot x = (P \cap G) \cdot x$.  This is because $P \cdot x$
splits into at most two $(P \cap G)$-orbits, and if it splits into two
orbits, then they must have the same dimension, namely $\dim \u$,
which is not possible.  Therefore, the representatives of Richardson
orbits given in Section \ref{S:orth} are valid for both $\OR_N(k)$
and $\SO_N(k)$.

Let $P$ be a parabolic subgroup of $G$ that is not conjugate to
$P(d) \cap G$ for any $d$.  Let $g \in \OR_N(k) \setminus G$, be
defined by $g = \sum_{i=2}^{n} (e_{i,i} + e_{-i,-i}) + e_{1,-1} +
e_{-1,1}$.  Then $gPg^{-1}$ is conjugate to $P(d) \cap \SO_N(k)$ for
some $d$, and can obtain a representative of the Richardson orbit
for $P$ from one for $gPg^{-1}$ by conjugating by $g$.
\end{rem}

\subsection{Symplectic groups}
\label{ss:symp} Much of the notation given in this subsection is
analogous to that for orthogonal groups in the previous subsection,
so we are briefer.

Let $G = \Sp_{2n}(k)$ and  $\g = \Lie G = \sp_{2n}(k)$. Let $V =
k^{2n}$ be the natural $2n$-dimensional $G$-module with standard
basis $v_{-n},\dots,v_{-1},v_1,\dots,v_n$ and $G$-invariant
skew-symmetric bilinear form $(\,,)$ defined by $(v_i, v_j) =
(v_{-i},v_{-j}) = 0$ and $(v_i, v_{-j}) = \delta_{i,j}$ for $1 \leq
i,j \leq n$.  Let $T$ be the maximal torus of $G$, which acts
diagonally on the standard basis of $V$ and let $\Phi$ be the root
system of $G$ with respect to $T$.  The following matrices give a
Chevalley basis for $\g$:
\begin{align*}
\{e_{i,j} - e_{-j,-i}\}_{1 \leq i,j \leq n} \cup
\{e_{i,-j}+e_{j,-i}, e_{-i,j}+e_{-j,i}\}_{1 \leq i<j \leq n} \cup
\{e_{k,-k}, e_{-k,k}\}_{1 \leq k \leq n},
\end{align*}
where $e_{i,j}$ denotes the $ij$-matrix unit. We define
$\epsilon_{i,j}$ to be the coefficient of $e_{i,j}$ in this basis if
it appears.

Let $d' = (d_1,\dots,d_s) \in \Z_{\ge 0}^s$ satisfy $\sum_{i=1}^s
d_i \le n$.  We define $d$ as in \S \ref{ss:orth}.  Such $d$ is
called a {\em symplectic dimension vector} and is said to be {\em
proper} if $d_i \ne 0$ for all $i \ne 0$. We define $c_i$ and $V_i$
as in \S \ref{ss:orth}. The stabilizer $P(d)$ of the isotropic flag
$0 = V_{-s} \sub V_{-s+1} \sub \dots \sub V_{-1}$ is a parabolic
subgroup of $G$. Moreover, for any parabolic subgroup $P$ of $G$,
there is a unique proper symplectic dimension vector
 $d$ such that $P$ is conjugate to $P(d)$.

In analogy with the orthogonal case, one can define a parabolic
subgroup $Q(d)$ of $\GL_N(k)$ such that $P(d) = Q(d) \cap G$.

We recall that the nilpotent $G$-orbits are parametrized by
partitions $\lambda$ of $N$ such that every odd part of $\lambda$
appears with even multiplicity. We call such a partition a {\em
symplectic partition}. Let $x$ and $y$ lie in the nilpotent
$G$-orbits parametrized by the symplectic partitions $\lambda$ and
$\mu$ respectively. Then $G \cdot x \sub \overline{G \cdot y}$ if
and only $\lambda \le \mu$, see \cite[\S 3]{H2}.

\section{Orthogonal groups} \label{S:orth}

For this section we use the notation given in \S \ref{ss:orth}. We
construct Richardson elements for all parabolic subgroups $P(d)$ of
$G = \OR_N(k)$. First we consider the case where all $d_i$ are at
most 2. For the general case we explain how to decompose $d$ as $d =
d^0 + d^1 + \dots + d^m$, where the entries of each $d^j$ are all 2
or less, in such a way that we can build up a Richardson element $x$
for $P(d)$ from the Richardson elements for the $P(d^j)$.  The idea
is that the natural $G$-module $V$ decomposes as an $x$-stable
orthogonal sum $V = V_1 \oplus \dots \oplus V_m$, where $P(d^j)$ is
a parabolic subgroup of $\OR(V_j)$.

\subsection{Blocks of size two or less} \label{ss:orth2}

Let $d = (d_{-s},\dots,d_s)$
be an
orthogonal dimension vector with $|d| = N$ and all $d_i = 0,1,2$,
further assume that if $d_0 = 1$, then all other nonzero entries of
$d$ are 1.  Let $P = P(d)$ be the parabolic subgroup of $G$
corresponding to $d$. In this subsection we construct a
representative $x \in \u$ of the Richardson orbit of $P$.

The idea is to define $x$ from a {\em line diagram} $D(d)$ in the
plane, which is defined by considering three cases.  The diagram
$D(d)$ consists of $N$ vertices labelled by the integers
$-n,-n+1,\dots,n-1,n = \lfloor \frac{N}{2} \rfloor$ (we omit 0 if
$N$ is even) and arrows between certain vertices. The vertices in
$D(d)$ are always labelled so that the labels increase from left to
right and from bottom to top; in order to save space we write $\bar
i$ rather than $-i$ in all the diagrams that we include below. We
define $x \in \u$ from $D(d)$ by
\begin{equation} \label{e:x}
x = \sum \epsilon_{i,j} e_{i,j},
\end{equation}
where the sum is taken over all arrows from $i$ to $j$ in $D(d)$ and
$\epsilon_{i,j}= \pm 1$ is as defined in \S \ref{ss:orth}.

The definition of $x$ means that $x$ sends $v_i$ to $\sum
\epsilon_{i,j} v_j$ where the sum is over all $j$ such that there is
an arrow from $i$ to $j$ in $D(d)$; this sum always has at most two
summands.  The diagram is centrally symmetric about the origin and
the labelling of the vertices forces the vertex labelled $-i$ to be
centrally symmetric to the vertex labelled $i$. This central
symmetry means that $x \in \g = \so_N$; the numbering of the
vertices ensures that $x \in \u$.

In each of the cases considered below, we give the partition
determined by the Jordan normal form of $x$; this is required for
the proof of Theorem \ref{T:orth} in \S \ref{ss:orthgen}.  As
notation for this we let $\sigma$ be the number of entries of $d$
equal to 2 and we let $\rho = |d| - \sigma$, so $\rho$ is the number
of nonzero entries in $d$.

\bigskip

{\em Case 1.} $d_0 = 1$ and all nonzero entries of $d$ are 1.

We construct the diagram $D(d)$ in the plane by placing vertices at
the points $(i,0)$ for $i$ such that $d_i = 1$. For each but
the leftmost vertex we draw an arrow to the vertex on its left.

For $d = (1,1,1,1,1,1,1,1,1)$ the diagram $D(d)$ is illustrated
below.
$$
\xymatrix{\bar 4 & \bar 3 \ar[l] & \bar 2 \ar[l] & \bar 1 \ar[l] & 0
\ar[l] & 1 \ar[l] & 2 \ar[l] & 3 \ar[l] & 4 \ar[l]}
$$
For $d = (1,0,1,0,1,0,1,0,1)$, the diagram $D(d)$ is as below.
$$
\xymatrix{\bar 2 & & \bar 1 \ar[ll] & & 0 \ar[ll] & & 1 \ar[ll] & &
4 \ar[ll]}
$$

We have $\rho = |d|$, $\sigma = 0$ and it is easy to see that the
Jordan normal form of $x$ is given by the partition $(\rho)$.

\bigskip

{\em Case 2.} $d_0 = 2$.

To construct $D(d)$, we draw vertices at $(i,1)$ and $(-i,-1)$ for $i
\ge 0$ such that $d_i \ne 0$; and at $(i,-1)$ and $(-i,1)$ for $i >
0$ such that $d_i =2$. From each vertex that is not leftmost in its
row we draw an arrow to the vertex on its left.  If some $d_i = 1$,
then we let $l > 0$ be minimal such that $d_l \ne 0$ and draw
additional arrows from $(l,1)$ to $(0,-1)$, and from $(0,1)$ to
$(-l,-1)$.

For example for $d = (2,2,2,2,2,2,2,2,2)$, we have $D(d)$ as below;
we recall that we label the vertices so that the labels increase
from left to right and bottom to top.
$$
\xymatrix{ \bar 8  & \bar 6 \ar[l] & \bar 4 \ar[l] & \bar 2 \ar[l] &
1 \ar[l] & 3 \ar[l] & 5 \ar[l] & 7 \ar[l] & 9 \ar[l]
\\ \bar 9 & \bar 7 \ar[l] & \bar 5 \ar[l] & \bar 3
\ar[l] & \bar 1 \ar[l] & 2 \ar[l] & 4 \ar[l] & 6 \ar[l] & 8 \ar[l]}
$$
We illustrate the diagram $D(d)$ when some entries of $d$ are 1 in
three examples.  For $d = (1,1,1,2,2,2,1,1,1)$ we have
$$
\xymatrix{& & & \bar 2 & 1 \ar[l] \ar[dl] & 3 \ar[l] \ar[dl] & 4
\ar[l] & 5 \ar[l] & 6 \ar[l] \\
\bar 6  & \bar 5 \ar[l] & \bar 4 \ar[l] & \bar 3 \ar[l] & \bar 1
\ar[l] & 2 \ar[l] }
$$
and for $d = (2,1,2,1,2,1,2,1,2)$ we have
$$
\xymatrix{\bar 6  & & \bar 3 \ar[ll] &  & 1 \ar[ll] \ar[dl] & 2
\ar[dl]
\ar[l] & 4 \ar[l] & 5 \ar[l] & 7 \ar[l] \\
\bar 7 & \bar 5 \ar[l] & \bar 4 \ar[l] & \bar 2 \ar[l] & \bar 1
\ar[l] & & 3 \ar[ll] & & 6 \ar[ll] }
$$
Finally we give an example when some entries of $d$ are 0, for $d =
(2,1,0,1,0,2,0,1,0,1,2)$ we have $l = 2$ and $D(d)$ as below.
$$
\xymatrix{\bar 4  & & & & & 1 \ar[lllll] \ar[dll] & & 2
\ar[ll] \ar[dll] & 3 \ar[l] & & 5 \ar[ll] \\
\bar 5 & \bar 3 \ar[l] & & \bar 2 \ar[ll] & & \bar 1 \ar[ll] & & & &
& 4 \ar[lllll]}
$$

Next we show that the Jordan normal form of $x$ is given by the
partition is $(\rho,\sigma)$.  In case all nonzero entries of $d$ are 2, we
have $\rho = \sigma = n$ and this is trivial.  So suppose, some
entries of $d$ are 1.  First, we note that the kernel of $x$ is
2-dimensional: if $d_0$ is the only entry of $d$ equal to 2 then the
kernel of $x$ has basis $\{e_{-n},e_1-e_{-1}\}$; otherwise, let $-a$
the label of the leftmost vertex on the top row of $D(d)$, then the
$\ker x$ has basis $\{e_{-n},e_{-a}\}$.  Therefore, $x$ has two
Jordan blocks. Next we note that $x^{\rho - 1} e_n = \pm 2e_{-n} \ne
0$ and $x^{\rho} = 0$ so that $x$ has a Jordan block of size $\rho$.
It follows that the partition given by the Jordan normal form of $x$
is $(\rho,\sigma)$ as required.

\bigskip

{\em Case 3.} $d_0 = 0$.

To construct $D(d)$, we draw vertices at $(i,1)$ and $(-i,-1)$ for $i
> 0$ such that $d_i \ne 0$; and at $(i,-1)$ and $(-i,1)$ for
$i > 0$ such that $d_i =2$.  From each vertex that is not leftmost
in its row we draw an arrow to the vertex on its left.  If four or
more entries of $d$ are 1, then there are additional diagonal arrows:
Let $0<l < m$ be minimal such that $d_l, d_m \ne 0$; then we draw an
arrow from $(l,1)$ to $(-m,-1)$, and an arrow from $(m,1)$ to
$(-l,-1)$.

If $d = (1,1,1,1,0,1,1,1,1)$, then $D(d)$ is the following diagram
$$
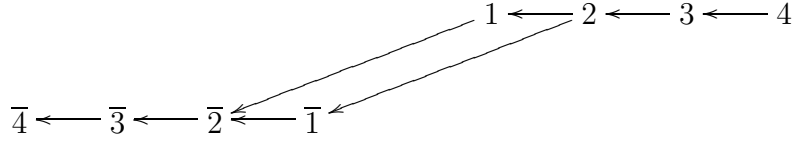
\xymatrix{& & & & & 1  \ar[dlll] & 2 \ar[l] \ar[dlll] & 3 \ar[l]
& 4 \ar[l] \\
\bar 4 & \bar 3 \ar[l] & \bar 2 \ar[l] & \bar 1 \ar[l]}
$$
and if $d = (2,2,2,2,0,2,2,2,2)$, then $D(d)$ is the following
diagram
$$
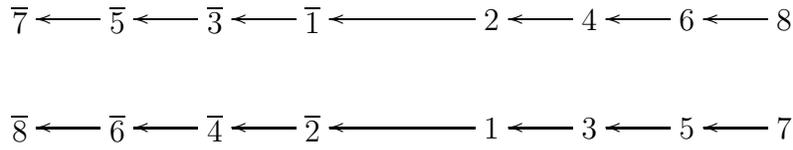
\xymatrix{ \bar 7 & \bar 5 \ar[l] & \bar 3 \ar[l] & \bar 1 \ar[l] &
& 2 \ar[ll] & 4 \ar[l] & 6 \ar[l] & 8 \ar[l]
\\ \bar 8 & \bar 6 \ar[l] & \bar 4 \ar[l] & \bar 2 \ar[l] & & 1 \ar[ll] & 3 \ar[l] & 5 \ar[l] &
7 \ar[l]}
$$
We illustrate the diagram $D(d)$ when some entries of $d$ are 1 in
the following examples. For $(2,2,2,1,2,2,1,2,2,2)$ we have
$$
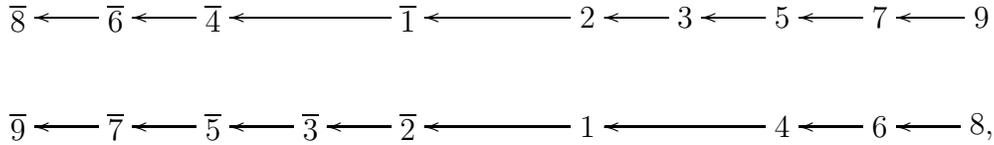
\xymatrix{\bar 8  & \bar 6 \ar[l] & \bar 4 \ar[l] & &  \bar 1
\ar[ll] & & 2 \ar[ll]
& 3 \ar[l] & 5 \ar[l] & 7 \ar[l] & 9 \ar[l] \\
\bar 9 & \bar 7 \ar[l] & \bar 5 \ar[l] & \bar 3 \ar[l]  & \bar 2
\ar[l] & & 1 \ar[ll] & & 4 \ar[ll] & 6 \ar[l] & 8, \ar[l] }
$$
for $d = (1,1,1,2,2,2,2,1,1,1)$ we have
$$
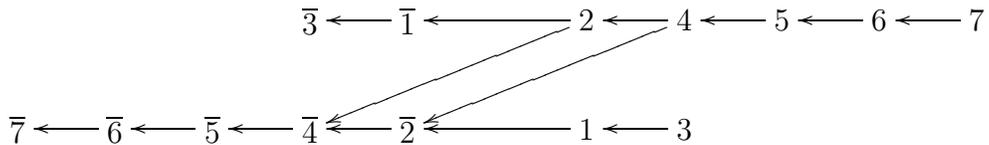
\xymatrix{& & & \bar 3 & \bar 1 \ar[l] & & 2 \ar[dlll] \ar[ll] & 4
\ar[dlll] \ar[l] & 5
\ar[l] & 6 \ar[l] & 7  \ar[l] \\
\bar 7  & \bar 6 \ar[l] & \bar 5 \ar[l] & \bar 4 \ar[l] & \bar 2
\ar[l] & & 1 \ar[ll] & 3 \ar[l] }
$$
and for $(2,1,2,1,2,2,1,2,1,2)$ we have
$$
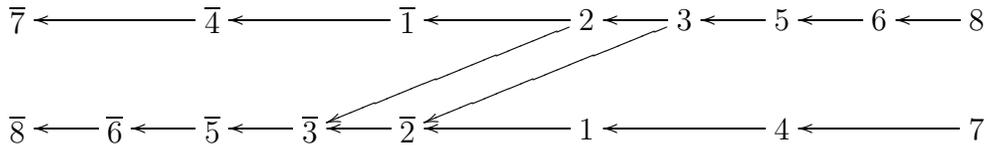
\xymatrix{\bar 7  & & \bar 4 \ar[ll] &  & \bar 1 \ar[ll] & & 2
\ar[dlll] \ar[ll]
& 3 \ar[dlll] \ar[l] & 5 \ar[l] & 6 \ar[l] & 8 \ar[l] \\
\bar 8 & \bar 6 \ar[l] & \bar 5 \ar[l] & \bar 3 \ar[l]  & \bar 2
\ar[l] & & 1 \ar[ll]  & & 4 \ar[ll] & & 7 \ar[ll] }
$$
Finally we give an example where some entries of $d$ are 0.  For $d
= (2,1,0,2,0,2,0,1,2)$ we have $l = 1$, $m=3$ and $D(d)$ as below.
$$
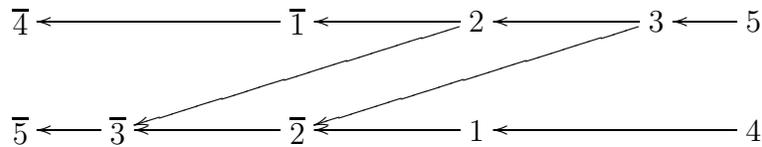
\xymatrix{\bar 4 & & & \bar 1 \ar[lll] & & 2
\ar[ll] \ar[dllll] & & 3 \ar[ll] \ar[dllll] & 5 \ar[l] \\
\bar 5 & \bar 3 \ar[l] & & \bar 2 \ar[ll] & & 1 \ar[ll] & & & 4
\ar[lll]}
$$

If all nonzero entries of $d$ are 2, then $\sigma = \rho$, and the
Jordan normal form of $x$ is given by the partition $(\rho,\rho)$.
Next we show that, when $d$ has some entries 1, the Jordan normal
form of $x$ is given by the partition is $(\rho-1,\sigma+1)$.  If
there are only two entries 1 in $d$ this is immediate ($\rho - 1 =
\sigma + 1$), so we assume that four or more entries of $d$ are 1.
First, we note that the kernel of $x$ is 2-dimensional: if all
nonzero entries of $d$ are 1, then the kernel of $x$ has basis
$\{e_{-n},e_1-e_{-1}\}$; otherwise, let $-a$ be the label of the
leftmost vertex on the top row of $D(d)$, then the $\ker x$ has
basis $\{e_{-n},e_{-a}\}$.  Next we note that $x^{\rho-2} e_n = \pm
2e_{-n} \ne 0$ and $x^{\rho-1} = 0$ so that $x$ has a Jordan block
of size $\rho-1$. It follows that the partition given by the Jordan
normal form of $x$ is $(\rho-1,\sigma+1)$ as required.

\subsection{General case} \label{ss:orthgen}

In this section we show how to build up Richardson elements for
arbitrary parabolic subgroups $P(d)$ of $G$.  This construction is
given in Definition \ref{D:orth} which is in two parts: first the
decomposition of the (proper) orthogonal dimension vector $d$ is given;
the second part gives the diagram
$D(d)$ from which $x$ is defined.  The Richardson element $x$ is
defined from the diagram $D(d)$ in the same way as in the previous
subsection from the formula \eqref{e:x}. The construction is quite
technical, so the reader may wish to look at the examples given
after the definition before reading it.

\begin{defn} \label{D:orth}
(i) First we explain how to make the decomposition $d = d^0 + d^1 +
\dots + d^m$.

We initially set $d^0 = 0$.  Suppose we have defined
$d^0,d^1,\dots,d^{j-1}$  and consider $c^j = d - \sum_{i=0}^{j-1}
d^i$. We define $d^j$ by considering cases:
\begin{enumerate}
\item[{\em Case A}.] $c^j_0 > 0$ and all nonzero entries of $c^j$ are 2 or greater.
Then we define $d^j$ by
$$
d^j_i = \left\{\begin{array}{cl} 0 & \text{if } c^j_i = 0 \\
                                 2 & \text{if } c^j_i \ge 2 \\
                \end{array}\right.
$$
\item[{\em Case B}.] $c^j_0 > 0$ is odd and $c^j$ has an entry equal to 1.
In this case we must have $d^0 = 0$ and we redefine it by
$$
d^0_i = \left\{\begin{array}{cl} 0 & \text{if } c^j_i = 0 \\
                                 1 & \text{if } c^j_i \ge 1 \\
                \end{array}\right.
$$
Then we update $c^j$.
\item[{\em Case C}.] $c^j_0 > 0$ is even and $c^j$ has an entry equal to 1.
We let $a \ge 0$
be the least positive even entry of $c^j$ and define $d^j$ by
$$
d^j_i = \left\{\begin{array}{cl} 0 & \text{if } c^j_i = 0 \\
                                 1 & \text{if } c^j_i
                                     \in\{1,3,\dots,a-1\} \\
                                 2 & \text{if } c^j_i \ge a \\
                \end{array}\right.
$$
\item[{\em Case D}.] $c^j_0 = 0$.  Then we define $d^j$ by
$$
d^j_i = \left\{\begin{array}{cl} 0 & \text{if } c^j_i = 0 \\
                                 1 & \text{if } c^j_i = 1 \\
                                 2 & \text{if } c^j_i \ge 2 \\
                \end{array}\right.
$$
\end{enumerate}
We continue until $c^j = 0$.

(ii) Now we explain how the diagram $D(d)$ is constructed. If $d_0$
is odd, then we draw the diagram $D(d^0)$ (with unlabelled
vertices); if $d_0$ is even then we must have $d^0 = 0$ and we do
nothing. Suppose we have dealt with $d^0,d^1,\dots,d^{j-1}$. Then we
insert the diagram $D(d^j)$ (with unlabelled vertices) stretched in
the vertical direction, so that a vertex at $(i,\pm 1)$ is moved to
$(i,\pm j)$ and diagonal arrows are stretched accordingly.

Label $D(d)$ so that numbers are increasing from left to right and
from bottom to top.  Then define $x$ from the formula \eqref{e:x}.
\end{defn}

\begin{exmp}
(1) Let $d = (3,4,2,4,3)$.  Then $d^0 = 0$, $d^1 = (2,2,2,2,2)$ and
$d^2 = (1,2,0,2,1)$.  The diagram $D(d)$ is as illustrated below.
$$
\xymatrix{& \bar 2 & & 5  \ar[ll] & 8 \ar[l] \\
          \bar 6 & \bar 3 \ar[l] & 1 \ar[l] & 4 \ar[l] & 7 \ar[l] \\
          \\
          \bar 7 & \bar 4 \ar[l] & \bar 1 \ar[l] & 3 \ar[l] & 6 \ar[l] \\
          \bar 8 & \bar 5 \ar[l] & & 2 \ar[ll] & }
$$

(2)  Let $d = (2,5,2,3,2,5,2)$.  Then we have $d^0 =
(0,1,0,1,0,1,0)$, $d^1 = (2,2,2,2,2,2,2)$ and $d^2 =
(0,2,0,0,0,2,0)$.  The diagram $D(d)$ is as illustrated below.
$$
\xymatrix{ & \bar 4 & & &  & 8 \ar[llll] \\
           \bar 9 & \bar 5 \ar[l] & \bar 2 \ar[l] & 1 \ar[l] & 3 \ar[l]  & 7 \ar[l]  & 10 \ar[l] \\
                  & \bar 6  & & 0 \ar[ll] & & 6 \ar[ll] & \\
           \bar 10 & \bar 7 \ar[l] & \bar 3 \ar[l] & \bar 1 \ar[l] & 2 \ar[l] & 5 \ar[l] & 9 \ar[l]\\
           & \bar 8 &  & &  & 4 \ar[llll] }
$$

(3)  Let $d = (2,2,4,1,4,2,2)$.  Then we have $d^0=(1,1,1,1,1,1,1)$,
$d^1 = (1,1,2,0,2,1,1)$ and $d^2 = (0,0,1,0,1,0,0)$.  The diagram
$D(d)$ is as illustrated below.
$$
\xymatrix{ & & & & 4 & \\
          & & \bar 1 & & 3 \ar[ll]\ar[ddlll] & 6 \ar[l]\ar[ddlll] &  8 \ar[l]  \\
          \bar 7&\bar 5\ar[l]& \bar 2\ar[l]& 0 \ar[l]& 2 \ar[l]& 5 \ar[l]& 7  \ar[l] \\
          \bar 8&\bar 6\ar[l]&\bar 3 \ar[l]& \ar[l]& 1 \ar[ll]&  \\
           & & \bar 4
}
$$

(4) Let $d = (4,1,3,4,3,1,4)$.  Then we have $d^0 = 0$, $d^1 =
(2,1,1,2,1,1,2)$ and $d^2 = (2,0,2,2,2,0,2)$.  The diagram $D(d)$ is
as illustrated below.
$$
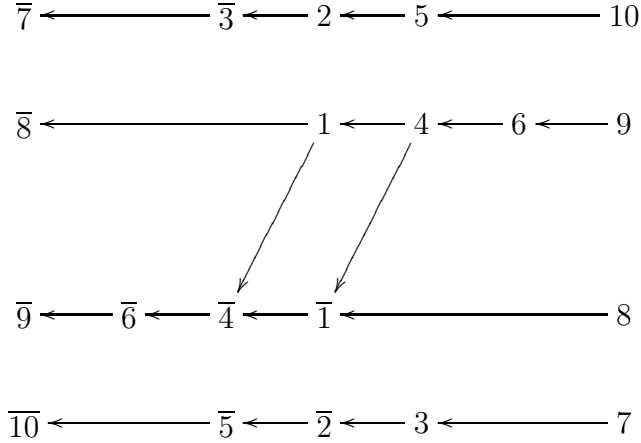
\xymatrix{\bar 7 & & \bar 3 \ar[ll] & 2 \ar[l] & 5 \ar[l] & & 10
\ar[ll] \\
\bar 8 & & & 1 \ar[lll] \ar[ddl] & 4 \ar[l] \ar[ddl] & 6\ar[l] & 9 \ar[l] \\
\\
\bar 9 & \bar 6 \ar[l] & \bar 4 \ar[l] & \bar 1 \ar[l] & & & 8 \ar[lll]\\
\bar{10} & & \bar 5 \ar[ll] & \bar 2 \ar[l] & 3 \ar[l] & & 7
\ar[ll]}
$$
\end{exmp}

Next we prove that our constructions do indeed give Richardson
elements.

\begin{thm} \label{T:orth}
The element $x$ defined from $D(d)$ by \eqref{e:x} is a Richardson
element for $P(d)$.
\end{thm}

\begin{proof}
Let $P = P(d)$ and let $Q = Q(d)$ be the corresponding parabolic
subgroup of $\GL_N(k)$, so we have $Q \cap \OR_N(k) = P$, see \S
\ref{ss:orth}. Let $\lambda = (\lambda_1,\lambda_2,\dots)$ be the
orthogonal partition of $N$ corresponding to the $\OR_N(k)$-orbit of
$x$ and let $\lambda^0 = (\lambda^0_1,\lambda^0_2\dots)$ be the
partition of $N$ corresponding to the Richardson orbit for $Q$. By
the construction of a Richardson element for $Q(d)$ given in
\cite{BHRR}, we know that $\lambda^0$ is the dual of the partition
of $N$; it is given by reordering the entries of $d$ into
non-increasing order.

We show that for any orthogonal partition $\lambda'$ of $N$
satisfying $\lambda' \le \lambda^0$ we have $\lambda' \le \lambda$
(in the dominance order on partitions). The closure order on nilpotent
$\GL_n(k)$-orbits and $G$-orbits (see \S \ref{ss:orth}) then tells
us that $\u \sub \overline{G \cdot x}$. Hence we have $(G \cdot x)
\cap \u$ is open in $\u$ so that $x$ is Richardson as required.

\smallskip

Let $a$ be the maximal such that $\lambda^0_1,\dots,\lambda^0_a$ are
all odd.  In case $N$ is odd we necessarily have $a$ is odd, and if
$N$ is even then so is $a$. We let $b = \lfloor \frac{a}{2}
\rfloor$.  From the construction of $D(d)$ in Definition
\ref{D:orth} and the partitions given for the $D(d^j)$ in Cases 1
and 2 in \S \ref{ss:orth2}, one can easily check that $\lambda_j =
\lambda^0_j$ for $j = 1,\dots,a$.

We define the partitions $\mu = (\lambda_{a+1},\lambda_{a+2},\dots)$
and $\mu^0 = (\lambda^0_{a+1},\lambda^0_{a+2},\dots)$. For $j =
1,\dots,m-b$, let $\sigma^j$ be the number of entries $2$ in
$d^{b+j}$ and let $\rho^j = |d^{b+j}| - \sigma^j$.  We define
$\delta^j= 0$ if all nonzero entries of $d^{b+j}$ are 2; and
$\delta^j = 1$ otherwise.  Then $\mu^0 =
(\rho^1,\sigma^1,\rho^2,\sigma^2,\dots)$ and $\mu =
(\rho^1-\delta^1,\sigma^1+\delta^1,\rho^2-\delta^2,\sigma^2+\delta^2,\dots)$.

We define a sequence $\lambda^0,\lambda^1,\dots,\lambda^{m-b} =
\lambda$ of orthogonal partitions of $N$ as follows: we define
$\lambda^j$ to have the same first $a+2j$ entries as $\lambda$ and
remaining entries equal to those of $\lambda^0$.  Below we show
inductively that for each $j$, we have that:
\begin{quote}
\begin{tabular}{p{350pt}l} $\lambda^j \le \lambda^0$ and for any orthogonal
partition $\lambda'$ of $N$ satisfying $\lambda' \le \lambda^0$ we
have $\lambda' \le \lambda^j$. & \raise-6pt\hbox{($\star$)}
\end{tabular}
\end{quote}
This condition for $j = m-b$ implies that $x$ is Richardson for $P$.

Assume inductively that ($\star$) holds for $\lambda^{j-1}$.  Let
$\lambda'$ be an orthogonal partition with $\lambda' \le \lambda$ in
the dominance ordering.  By induction we have that $\sum_{i = 1}^l
\lambda'_i \le \sum_{i = 1}^l \lambda^j_i$ for all $l \ne a+2j-1$.
The only way that we can have $\sum_{i = 1}^{a+2j-1} \lambda'_i >
\sum_{i = 1}^{a+2j-1} \lambda^j_i$ is if $\delta^j = 1$, $\lambda'_i
= \lambda^j_i$ for $i = 1,\dots,a+2j-2$ and $\lambda'_{a+2j-1} =
\lambda^0_{a+2j-1} = \lambda^j_{a+2j-1} + 1$.  However, these
conditions force $\lambda^0_{a+2j-1}$ to be even and of odd
multiplicity in $\lambda'$, contrary the assumption that $\lambda'$
is an orthogonal partition. It follows that ($\star$) holds for
$\lambda^j$.
\end{proof}

As mentioned after the statement of Theorem \ref{T:support} we now
discuss the type of the Richardson element $x$.
Let $\Gamma$ be the subset of $\Phi$
such that $x = \sum_{\alpha \in \Gamma} e_\alpha$; the elements of
$\Gamma$ correspond to pairs of arrows in $D(d)$.

First we consider the case when $N$ is even.  For $j = 1,\dots,m$,
we let $\sigma^j$ be the number of entries $2$ in $d^j$ and let
$\rho^j = |d^j| - \sigma^j$.  We define $\delta^j= 0$ if $d^j_0 = 2$
or all nonzero entries of $d6j$ are 2; and $\delta^j = 1$ otherwise.
We let $I_d = \{j \mid \rho^j - \delta^j = \sigma^j + \delta^j\}$,
$J_d = \{j \mid \rho^j - \delta^j > \sigma^j + \delta^j\}$ and $\eta
= \frac{1}{2}\sum_{j \in J_d} (\rho^j + \sigma^j)$. Then the type of
the orbit of $x$ is
$$
D_\eta + \sum_{j \in I_d} A_{\rho^j - \delta^j - 1}.
$$
By convention we have $D_2 = A_1 + A_1$  (respectively $D_3 = A_3$),
but the corresponding Levi subgroup lies in a different conjugacy
class to an $A_1 + A_1$ (respectively $A_3$) Levi subgroup. Further,
we note that if $\eta = 0$ and $|I_d| = n - \sum_{j \in I_d} (\rho^j
- \delta^j - 1)$, then inside $\SO_N(k)$ there are two conjugacy classes
of Levi subgroups with isomorphism type $\sum_{j \in I_d} A_{\rho^j
- \delta^j}$; however, these two classes fuse in $G = \OR_N(k)$.
The type of $x$ given above can be easily verified as in \cite[\S
3]{Pa} or \cite[\S 4.8]{Ja}.

For $j \in I_d$, there are $2(\rho^j - \delta^j -1)$ arrows in
$D(d^j)$. Therefore, the arrows in $D(d^j)$ contribute $\rho^j -
\delta^j -1$ elements to $\Gamma$. For $j \in I_d$, there are
$\rho^j + \sigma^j$ arrows in $D(d^j)$.  Therefore, the arrows in
$D(d^j)$ contribute $\frac{1}{2}(\rho^j + \sigma^j)$ elements to
$\Gamma$. Hence, we have
$$
|\Gamma| = \eta + \sum_{j \in I_d} (\rho^j - \delta^j - 1).
$$
This is the minimal possible size for $\Gamma$ as stated in Theorem
\ref{T:support}, because it is the rank of the type of the orbit of
$x$.

\smallskip

The situation for the case where $N$ is odd is almost exactly the
same.  We define $\rho^j$, $\sigma^j$, $\delta^j$, $I_d$ and $J_d$
as for the case $N$ even.  We define $\eta = \frac{1}{2}(|d^0|-1 +
\sum_{j \in J_d} (\rho^j + \sigma^j))$. Then the type of the orbit
of $x$ is of the form
$$
B_\eta + \sum_{j \in I_d} A_{\rho^j - \delta^j - 1},
$$
where by convention $B_1 = A_1$, but with a short root and therefore
the corresponding Levi subgroup is in a different conjugacy class to
an $A_1$ Levi subgroup. One can check this type is correct, and that
the size of $\Gamma$ is equal to the rank of the type of $x$, as for
the case $N$ even.

\section{Symplectic groups} \label{S:symp}

In this section we use the notation from \S \ref{ss:symp}.  We
construct Richardson elements for all parabolic subgroups $P(d)$ of
$G= \Sp_{2n}(k)$. As for the orthogonal groups, we first consider
the case where all $d_i$ are at most 2.  We decompose a general
symplectic dimension vector $d$ as $d = d^1 + \dots + d^m$ where the
entries of each $d^j$ are all 2 or less, and then build up a
Richardson element for $P(d)$ from the Richardson elements for the
$P(d^j)$.

\subsection{Block of size two or less}

Let $d = (d_{-s},\dots,d_s)$
be a
symplectic dimension vector with $|d| = 2n$ and all $d_i = 0,1,2$.
Let $P = P(d)$ be the corresponding parabolic subgroup of $G$. In
this subsection we construct a representative $x \in \u$ of the
Richardson orbit of $P$.

As for the orthogonal case, we define $x$ from a {\em line diagram}
$D(d)$ in the plane which is given by considering four cases. The
diagram consists of vertices labelled $\pm 1,\pm 2,\dots,\pm n$ and
arrows between certain vertices. The vertices in $D(d)$ are labelled
as in the previous section: they increase from left to right and
from bottom to top; again $\bar i$ stands for $-i$. We define $x \in
\u$ by
\begin{equation}\label{e:x-sp}
x = \sum \epsilon_{i,j} e_{i,j},
\end{equation}
where the sum is taken over all arrows from $i$ to $j$ in $D(d)$ and
$\epsilon_{i,j} = \pm 1$ is defined as in \S \ref{ss:symp}.

The diagram $D(d)$ is centrally symmetric about the origin and the
vertex labelled $-i$ is centrally symmetric to the vertex $i$, which
implies $x\in\sp_{2n}$.  Thanks to the numbering of the vertices we
have $x\in\u$.

For each of the four cases below we give the partition of the Jordan
normal form of $x$. As notation for this we let $\sigma$ be the
number of entries of $d$ equal to 2 and we let $\rho = |d| -
\sigma$, so $\rho$ is the number of nonzero entries in $d$.

\bigskip

{\em Case 1.} $d_0=0$ and all nonzero entries of $d$ are 1.

To construct $D(d)$ we draw vertices at the points $(i,0)$ for $i$
with $d_i \ne 0$. From each vertex that has a left neighbour, we
draw an arrow to the left.

So for example, for $d = (1,1,1,1,0,1,1,1,1)$, we have $D(d)$ as
below
$$
\xymatrix{\bar 4 & \bar 3 \ar[l] & \bar 2 \ar[l] & \bar 1 \ar[l] & &
1 \ar[ll] & 2 \ar[l] & 3 \ar[l] & 4 \ar[l]}
$$
and for $d = (1,0,1,1,0,0,0,1,1,0,1)$ we have
$$
\xymatrix{\bar 3 & & \bar 2 \ar[ll] & \bar 1 \ar[l] & & & & 1
\ar[llll] & 2 \ar[l] & &  3 \ar[ll]}
$$

We have $\rho = |d|$ and $\sigma = 0$, and the Jordan normal form of
$x$ is given by the partition $(\rho)$.

\bigskip

{\em Case 2.} $d_0 = 2$.

We draw vertices at $(i,1)$ and $(-i,-1)$ for $i \ge 0$ such that
$d_i \ne 0$ and further vertices at $(-i,1)$ and at $(i,-1)$ for $i
> 0$ such that $d_i=2$. From each vertex that is not leftmost in its
row we draw an arrow to the vertex on its left.  If four or more
entries of $d$ are 1, then let $l > 0$ be minimal such that
$d_l > 0$ and draw an additional arrow  from
$(l,1)$ to $(-l,-1)$.

So if $d = (2,2,2,2,2,2,2)$, then we have $D(d)$ as below
$$
\xymatrix{\bar 6  & \bar 4 \ar[l] & \bar 2 \ar[l]  & 1 \ar[l]
& 3 \ar[l] & 5 \ar[l] & 7 \ar[l] \\
\bar 7 & \bar 5 \ar[l] & \bar 3 \ar[l] & \bar 1 \ar[l] & 2\ar[l] & 4
\ar[l] & 6 \ar[l]}
$$
We illustrate $D(d)$ in cases where there are some entries 1 in $d$
with three examples.  For (2,1,2,2,2,1,2) we have
$$
\xymatrix{\bar 5  &  &  \bar 2 \ar[ll] &  1 \ar[l] & 3 \ar[l] & 4
\ar[l]
 & 6 \ar[l]  \\
\bar 6 & \bar 4\ar[l] & \bar 3 \ar[l] & \bar 1\ar[l] & 2 \ar[l] & &
5 \ar[ll] }
$$
for (1,1,2,2,2,2,2,1,1) we have
$$
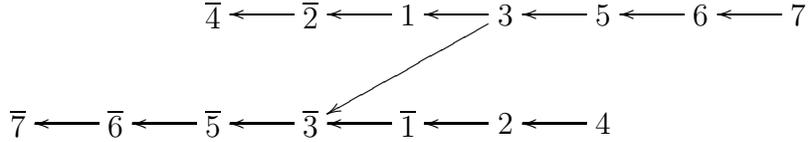
\xymatrix{& & \bar 4 & \bar 2 \ar[l] & 1 \ar[l] & 3
\ar[l]  \ar[dll] & 5 \ar[l] & 6 \ar[l] & 7 \ar[l]  \\
\bar 7  & \bar 6\ar[l] & \bar 5 \ar[l] & \bar 3 \ar[l] & \bar 1
\ar[l] & 2 \ar[l] & 4 \ar[l] }
$$
and for (2,1,1,2,1,1,2) we have
$$
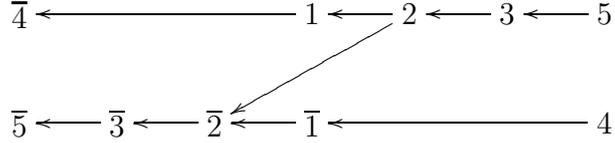
\xymatrix{  \bar 4  & & & 1 \ar[lll] & 2 \ar[l] \ar[dll] & 3
\ar[l] & 5 \ar[l] \\
\bar 5  & \bar 3 \ar[l] & \bar 2 \ar[l] & \bar 1 \ar[l] & & & 4
\ar[lll] }
$$
Finally, we give an example when some of the entries of $d$ are 0.
For $(1,0,2,1,0,2,0,1,2,0,1)$ we have $l = 2$ and $D(d)$ as below
$$
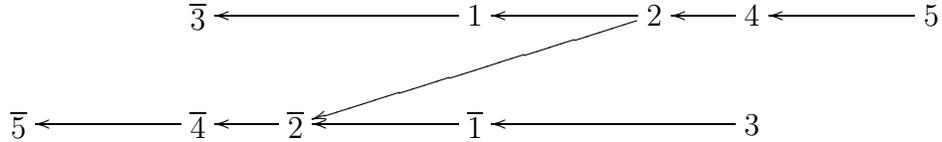
\xymatrix{ & & \bar 3 & & & 1 \ar[lll] & & 2 \ar[ll] \ar[dllll] & 4 \ar[l] & & 5 \ar[ll] \\
\bar 5 & & \bar 4 \ar[ll] & \bar 2 \ar[l] & & \bar 1 \ar[ll] & & & 3
\ar[lll]}
$$

If all entries of $d$ are 2 then $\rho = \sigma$ and the Jordan
normal form of $x$ is given by the partition $(\rho,\rho)$;
otherwise the Jordan normal form is given by $(\rho - 1,\sigma +1)$.
This can be verified using argument similar to those for orthogonal
cases in \S \ref{ss:orth2}.

\bigskip

{\em Case 3.}  $d_0 = 0$ and all nonzero entries of $d$ are 2.

To construct $D(d)$ we draw vertices at the points $(i,1)$ and
$(i,-1)$ for $i=-s,\dots,s$ such that $d_i = 2$.
From each vertex which is not
leftmost in its row we draw an arrow to the vertex on its left.

For example for $d = (2,2,2,0,2,2,2)$, we have $D(d)$ as below.
$$
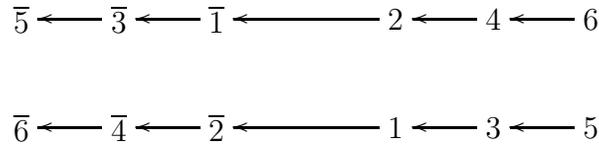
\xymatrix{\bar 5 & \bar 3 \ar[l] & \bar 1 \ar[l] & & 2\ar[ll] & 4 \ar[l] & 6 \ar[l] \\
\bar 6 & \bar 4 \ar[l] & \bar 2 \ar[l] & & 1 \ar[ll] & 3\ar[l] & 5
\ar[l]}
$$
For $d = (2,0,2,0,0,0,2,0,2)$ we have
$$
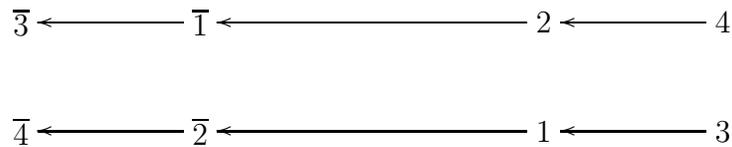
\xymatrix{\bar 3 & & \bar 1 \ar[ll] & & & & 2 \ar[llll] & & 4 \ar[ll] \\
\bar4 & & \bar 2 \ar[ll]  & & & & 1 \ar[llll] & & 3 \ar[ll]}
$$

In this case we have $\rho = \sigma$ and the partition given by the
Jordan normal form of $x$ is $(\rho,\rho)$.

\bigskip

{\em Case 4.}  $d_0=0$ and $d$ has both entries $1$ and $2$.

To construct $D(d)$ we draw vertices at the points $(i,1)$ and
$(-i,-1)$ for $i > 0$ such that $d_i \ne 0$; furthermore, for $i>0$
with $d_i=2$, we draw vertices at $(i,-1)$ and $(-i,1)$.  We let
$l,m > 0$ be minimal such that $d_l \ne 0$ and $d_m = 2$ (allowing
$l = m$). From each vertex that is not leftmost in its row and is
not at $(l,1)$ or $(m,-1)$, we draw an arrow to the vertex on its
left. We draw additional arrows from $(l,1)$ to $(-l,-1)$ and from
$(m,-1)$ to $(-m,1)$.

We illustrate $D(d)$ in three examples: For $d =
(1,1,2,2,0,2,2,1,1)$ we have $l = m = 1$ and $D(d)$ as below
$$
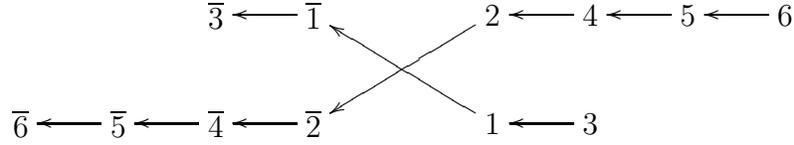
\xymatrix{& & \bar 3  & \bar 1 \ar[l]  & & 2 \ar[dll] & 4 \ar[l] & 5
\ar[l] & 6 \ar[l] \\
\bar 6  & \bar 5 \ar[l] & \bar 4 \ar[l] & \bar
2 \ar[l] & & 1 \ar[ull] & 3 \ar[l] }
$$
for $d = (2,2,1,1,0,1,1,2,2)$ we have $l=1$, $m = 3$ and $D(d)$ as
below
$$
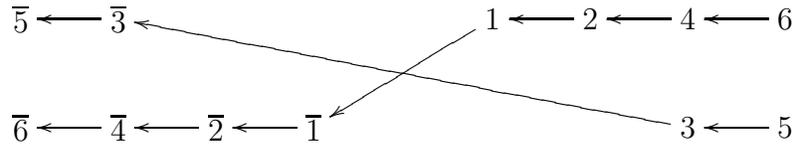
\xymatrix{\bar 5  & \bar 3 \ar[l] & & & & 1 \ar[dll] & 2
\ar[l] & 4 \ar[l] & 6 \ar[l] \\
\bar 6 & \bar 4 \ar[l] & \bar 2 \ar[l] & \bar 1 \ar[l] & & & & 3
\ar[ullllll] & 5  \ar[l] }
$$
and for $d = (1,0,1,2,1,0,0,0,1,2,1,0,1)$ we have $l=2$, $m=3$ and
$D(d)$ as below
$$
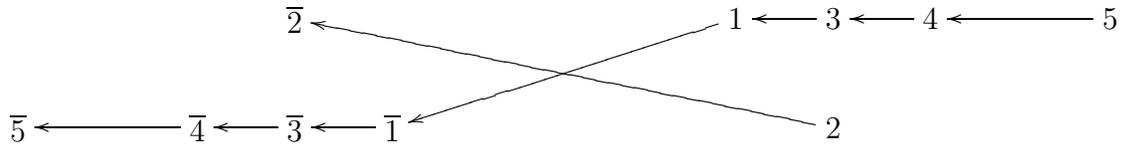
\xymatrix{ & & & \bar 2 & & & & & 1 \ar[dllll] & 3 \ar[l]& 4 \ar[l]& & 5 \ar[ll] \\
\bar 5 & & \bar 4 \ar[ll] & \bar 3 \ar[l] & \bar 1 \ar[l] & & & & &
2 \ar[ullllll]}
$$

It is easy to see that the Jordan normal form of $x$ is given by the
partition $(\rho,\sigma)$.

\subsection{General case} \label{ss:sympgen}
Analogously to the orthogonal groups, we can now build up Richardson
elements for arbitrary subgroups $P(d)$ of $G$.  Given a (proper)
symplectic dimension vector $d$, we decompose it as a sum of
dimension vectors $d^j$ with entries at most 2. Then we build up a
diagram $D(d)$ from the diagrams $D(d^j)$, which defines the
Richardson element $x$ as in the previous subsection from the
formula \eqref{e:x-sp}. After describing the details of the
construction we illustrate it with examples.

\begin{defn}
(i) First we explain how to make the decomposition $d = d^0 + d^1 +
\dots + d^m$.

We initially set $d^0 = 0$.  Suppose we have defined
$d^0,d^1,\dots,d^{j-1}$  and consider $c^j = d - \sum_{i=0}^{j-1}
d^i$. We define $d^j$ by considering cases:
\begin{enumerate}
\item[{\em Case A}.] $c^j_0 > 0$.
Then we define $d^j$ by
$$
d^j_i = \left\{\begin{array}{cl} 0 & \text{if } c^j_i = 0 \\
                                 1 & \text{if } c^j_i = 1 \\
                                 2 & \text{if } c^j_i \ge 2 \\
                \end{array}\right.
$$
\item[{\em Case B}.]  $c^j_0 = 0$ and all nonzero entries of $c^j$ are 2 or greater.
Then we define
$d^j$ by
$$
d^j_i = \left\{\begin{array}{cl} 0 & \text{if } c^j_i = 0 \\
                                 2 & \text{if } c^j_i \ge 2 \\
                \end{array}\right.
$$
\item[{\em Case C}.]  $c^j_0 = 0$ and $c^j$ has an entry equal to 1
and a positive even entry. Then we let $a > 0$ be the least positive
even entry of $c^j$ and define $d^j$ by
$$
d^j_i = \left\{\begin{array}{cl} 0 & \text{if } c^j_i = 0 \\
                                 1 & \text{if } c^j_i
                                     \in\{1,3,\dots,a-1\} \\
                                 2 & \text{if } c^j_i \ge a \\
                \end{array}\right.
$$
\item[{\em Case D}.]  $c^j_0 = 0$ and all nonzero entries of $c^j$
are odd.  In this case we must have $d^0 = 0$ and we redefine it by
$$
d^0_i = \left\{\begin{array}{cl} 0 & \text{if } c^j_i = 0 \\
                                 1 & \text{if } c^j_i \ge 1 \\
                \end{array}\right.
$$
Then we update $c^j$.
\end{enumerate}
We continue until $c^j = 0$.
%
%
%

(ii) Now the diagram $D(d)$ is constructed as follows. If $d_0 \ne
0$, then we draw the diagram $D(d^0)$ (with unlabelled vertices); if
$d^0 = 0$ then we do nothing. Suppose we have taken care of
$d^0,d^1,d^2,\dots,d^{j-1}$. Then we insert the diagram $D(d^j)$
(with unlabelled vertices) stretched in the vertical direction, so
that a vertex $(i,\pm 1)$ is moved to $(i,\pm j)$ and diagonal
arrows are stretched accordingly.

Label $D(d)$ with the numbers $\pm 1,\dots,\pm n$ increasing from
left to right and bottom to top. Then define $x$ from the formula
\eqref{e:x-sp}
\end{defn}

\begin{exmp}
(1) Let $d=(3,4,2,4,3)$. Then $d^0 = 0$, $d^1=(2,2,2,2,2)$ and
$d^2=(1,2,0,2,1)$. The diagram $D(d)$ is shown below.
$$
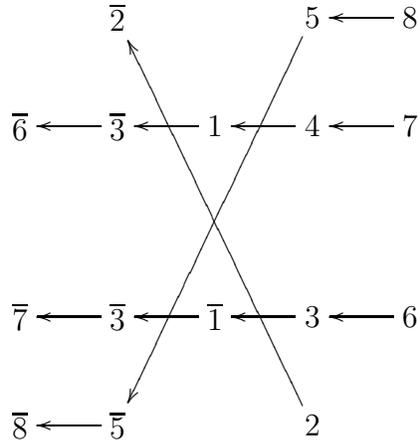
\xymatrix{   & \bar 2 & & 5\ar[lldddd] & 8\ar[l] \\
\bar 6  & \bar 3\ar[l] & 1\ar[l] & 4 \ar[l] & 7\ar[l] \\
 & & \\
\bar 7 & \bar 3\ar[l] & \bar 1\ar[l] & 3\ar[l] & 6\ar[l] \\
\bar 8  & \bar 5\ar[l]  & & 2\ar[lluuuu]  }
$$

(2)  Let $d = (2,3,2,2,2,3,2)$.  Then we have $d^0=(0,1,0,0,0,1,0)$
and $d^1=(2,2,2,2,2,2,2)$.  The diagram $D(d)$ is as illustrated
below.
$$
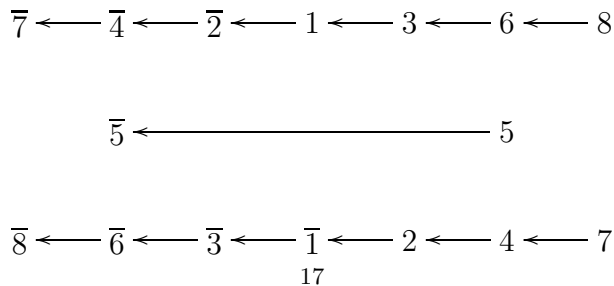
\xymatrix{\bar{7} & \bar 4 \ar[l] & \bar 2 \ar[l] & 1 \ar[l]&  3\ar[l] & 6\ar[l] & 8\ar[l] \\
                   & \bar 5 &  &  & & 5 \ar[llll] &  \\
          \bar{8} & \bar{6} \ar[l] & \bar 3 \ar[l] & \bar 1 \ar[l] & 2\ar[l] & 4 \ar[l]& 7 \ar[l] \\
          }
$$


(3)  Let $d=(3,1,6,1,1,2,1,1,6,1,3)$. Then $d^0 = 0 $,
$d^1=(2,1,2,1,1,2,1,1,2,1,2)$, $d^2=(1,0,2,0,0,0,0,0,2,0,1)$ and
$d_3 = (0,0,2,0,0,0,0,0,2,0,0)$. The diagram $D(d)$ is shown below.
$$
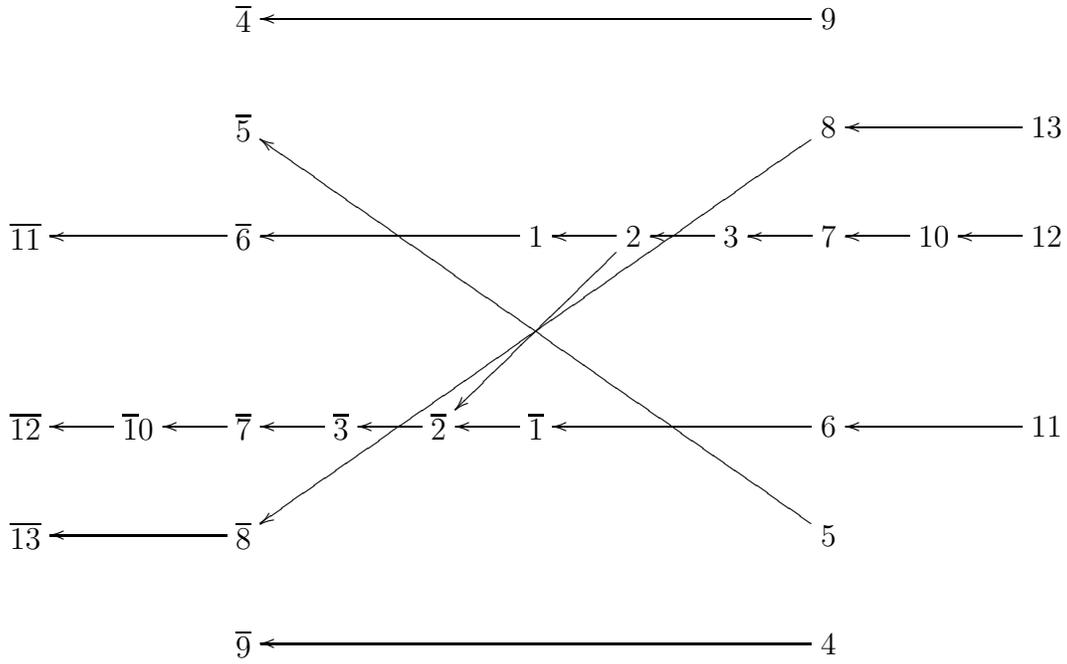
\xymatrix{ & & \bar 4 & & & & &  & 9 \ar[llllll] \\
           & & \bar 5 & & & & & & 8 \ar[ddddllllll] & & 13 \ar[ll]\\
           \bar{11} & & \bar 6 \ar[ll] & & & 1 \ar[lll]  & 2 \ar[l] \ar[ddll] & 3 \ar[l] & 7 \ar[l]  & 10 \ar[l] & 12 \ar[l] \\
           \\
           \bar{12} & \bar 10 \ar[l] & \bar 7 \ar[l] & \bar 3 \ar[l] & \bar 2 \ar[l] & \bar 1 \ar[l] & &  & 6 \ar[lll] & & 11 \ar[ll] \\
           \bar{13} & & \bar 8 \ar[ll] & & & & & & 5 \ar[uuuullllll] \\
         &  & \bar 9  &  & & & &  & 4 \ar[llllll] }
$$

(4)  Let $d = (5,3,5,0,5,3,5)$.  Then $d^0 = (1,1,1,0,1,1,1)$, $d^1
= (2,2,2,0,2,2,2)$ and $d^2 = (2,0,2,0,2,0,2)$.  the diagram $D(d)$
is shown below.
$$
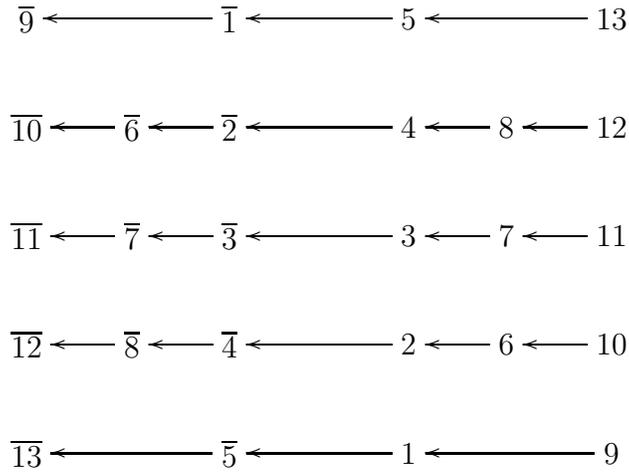
\xymatrix{ \bar 9 & & \bar 1 \ar[ll] & & 5 \ar[ll] & & 13 \ar[ll]
\\
\bar{10} & \bar 6 \ar[l] & \bar 2 \ar[l] & & 4 \ar[ll] & 8 \ar[l] &
12 \ar[l] \\
\bar{11} & \bar 7 \ar[l] & \bar 3 \ar[l] & & 3 \ar[ll] & 7 \ar[l] &
11 \ar[l] \\
\bar{12} & \bar 8 \ar[l] & \bar 4 \ar[l] & & 2 \ar[ll] & 6 \ar[l] &
10 \ar[l] \\
\bar{13} & & \bar 5 \ar[ll] & & 1 \ar[ll] & & 9 \ar[ll]
\\}
$$
\end{exmp}

We now state a theorem saying that the Richardson elements
constructed as in Definition \ref{D:orth} are indeed Richardson
elements. The proof of Theorem \ref{T:symp} is completely analogous
to the proof of the Theorem \ref{T:orth}; it is based on the fact
that the closure order for nilpotent orbits is given by the
dominance order on partitions as mentioned in \S \ref{ss:symp}.
Therefore, we omit the details.

\begin{thm} \label{T:symp}
The element $x$ defined from $D(d)$ is a Richardson element for
$P(d)$.
\end{thm}

As mentioned after the statement of Theorem \ref{T:support} we now
discuss the type of the Richardson element $x$.
Let $\Gamma$ be the subset of $\Phi$ such that $x = \sum_{\alpha \in
\Gamma} e_\alpha$; the elements of $\Gamma$ correspond to either
single arrows in $D(d)$ that pass through the origin, or pairs of
arrows in $D(d)$ that do not pass through the origin.  One can check
the type of $x$ given below is correct, and that it has rank equal
to the size of $\Gamma$ as for the orthogonal case in \S
\ref{ss:orthgen}.  Therefore, we omit the details.

For $j = 1,\dots,m$, we let $\sigma^j$ be the number of entries $2$
in $d^j$ and let $\rho^j = |d^j| - \sigma^j$.  We define $\delta^j=
0$ if $d^j_0 = 0$ or all nonzero entries of $d^j$ are 2; and $\delta^j
= 1$ otherwise.  We let $I_d = \{j \mid \rho^j - \delta^j = \sigma^j
+ \delta^j\}$, $J_d = \{j \mid \rho^j - 1
> \sigma^j - 1\}$ and $\eta = \frac{1}{2}(|d^0| + \sum_{j \in
J_d} (\rho^j + \sigma^j))$. Then the type of the orbit of $x$ is
$$
C_\eta + \sum_{j \in I_d} A_{\rho^j - \delta^j - 1},
$$
where by convention $C_1 = A_1$, but with a long root and therefore
the corresponding Levi subgroup is in a different conjugacy class to
an $A_1$ Levi subgroup.

\end{document}